\documentclass{article} 

\usepackage[english]{babel} 
\usepackage{amssymb}
\usepackage{amsmath}
\usepackage{txfonts}
\usepackage{mathdots}
\usepackage[classicReIm]{kpfonts}
\usepackage{graphicx}

\usepackage{tocloft}

\title{On the distribution of composite odd numbers\\}
\author{\textbf{WOLF Marc, WOLF Fran\c{c}ois, VILLEMIN Fran\c{c}ois-Xavier} \\ \\
Email: marc.wolf3@wanadoo.fr \\
Email: francois.wolf@dbmail.com \\
Email: fxvillemin@orange.fr \\
\\
\textbf{\underbar{August 26, 2018}}\\
}
\date{}

\begin{document}


\maketitle

\begin{abstract}
We study odd numbers through a straightforward indexing. We focus in particular on odd prime and composite numbers and their distribution. With a counting argument, we calculate the limit of two sums and compare their convergence rate.\\
\\
\textbf{\underbar{Keywords:}} composite odd numbers, prime numbers, M\"{o}bius function, prime number theorem, counting, alternating sums, convergence rate, Euler approximation.

\end{abstract}

\tableofcontents

\noindent 
\section*{Introduction}
\addcontentsline{toc}{part}{Introduction}
All prime numbers greater than 3 are of the form $6m-1$ or $6m\ +\ 1$. This condition is obviously not sufficient to characterize them, and we propose to study in this article the distribution of composite numbers among them. By counting arguments, we will calculate the limit of a particular alternating sum.

\noindent 
\section{Indices of composite odd numbers: the set $\boldsymbol{W}$}

\noindent Let us recall here some of the notations from [1].

\noindent 1. $I$ is the set of odd integers greater than 1, i.e.:
\[I=\left\{N_k=2k+3;k\in \mathbb{N}\right\}\] 
with $k$ the \underbar{index} of odd number $N_k$;

\noindent 2. $P$ is the set of odd prime numbers, primes will also be enumerated in ascending order $p_0=3,p_1,\dots ,p_n\dots $, with $q_0=0,q_1,\dots ,q_n\dots $ their respective indices;

\noindent 3. $C$ is the set of composite odd integers greater than 1, i.e.:
\[C=I\backslash P=\{N_k\in I|\exists \left(a,b\right)\in I,\ N_k=ab\}\] 
The function $f:k\in \mathbb{N}\mathrm{\longmapsto }N_k\in I$ is bijective. The inverse function is $f^{-1}:N_k\mathrm{\in }I\mathrm{\longmapsto }\mathrm{\ }\ k=\frac{N_k-3}{2}$. The inverse image of $C$ is called $W$:
\[W=f^{-1}\left(C\right)=\{k\in \mathbb{N}\mathrm{|}\ N_k\in C\}\] 
It is reminded that the function
\[k:\left(n,j\right)\in {\mathbb{N}}^*\mathrm{\times }\mathbb{N}\mapsto k_j\left(n\right)=\left(2j+3\right)n+j\] 
is a surjection on $W$. In other words, $W$ is the (non-disjoint) union of the sets $W_j={\left\{k_j\left(n\right)\right\}}_{n\in {\mathbb{N}}^*}$.

\noindent 

\noindent Finally, the remarkable indices are the indices of the form $k_j\left(j+1\right)=2j^2+6j+3$, they correspond to the indices of odd squares.

\noindent \underbar{}

\noindent 
\subsection{Partitioning the set of odd numbers}

\noindent Every odd prime number greater than 3 is of the form $6m-1$ or $6m+1$.

\noindent Thus, the indices of odd prime numbers greater than 3 must be respectively of the form $3n+1\ $and $3n+2$. This means that the index of an odd number $x$ is a multiple of 3 if and only if $x$ itself is as well.

\noindent 

\noindent \textbf{\underbar{Definitions 1.1}}\underbar{:}

\begin{enumerate}
\item \underbar{ }For all $j\in \mathbb{N}$, the interval $I_j$ is defined by:
\[I_j=\left\{ \begin{array}{c}
\llbracket 0,11 \rrbracket\ \mathrm{if\ }j=0\ \ \ \ \ \ \ \ \ \ \ \ \ \ \ \ \ \ \ \ \ \ \ \ \ \ \ \ \ \ \ \ \ \ \ \ \ \ \ \ \ \ \ \ \  \\ 
\llbracket k_j\left(j+1\right)+1,k_{j+1}\left(j+2\right)\rrbracket \ \mathrm{otherwise.} \end{array}
\right.\] 
\end{enumerate}

\begin{enumerate}
\item  Let the unit $U\left(j\right)$ be the size of $I_j$, i.e.:
\[U\left(0\right)=12\ \mathrm{and\ for\ all\ }j>0,\ U\left(j\right)=4j+8\] 
\end{enumerate}

\begin{enumerate}
\item  The counting interval $I_D\left(j\right)$ is:
\[I_D\left(j\right)\mathrm{=}\llbracket 0,k_{j+1}\left(j+2\right)\rrbracket =\bigcup^j_{q=0}{I_q}.\] 
\end{enumerate}
The size of a set $E$ is noted $\left|E\right|$.

\noindent The former union equal to $I_D\left(j\right)$ being clearly disjoint, we easily deduce the following equality:
\[\left|I_D\left(j\right)\right|=\sum^j_{q=0}{U\left(q\right)}=12+\sum^j_{q=1}{4q+8}=k_{j+1}\left(j+2\right)+1=2j^2+10j+12.\]

\noindent 
\section{Indices of the form $\boldsymbol{3}\boldsymbol{n}\boldsymbol{+}\boldsymbol{1}$ and $\boldsymbol{3}\boldsymbol{n}\boldsymbol{+}\boldsymbol{2}$}

\noindent 
\subsection{The set A of odd number indices and B of composite odd number indices}

\noindent \underbar{}

\noindent \textbf{\underbar{Definition 2.1.1}}: Let $A\left(j\right)$ be the set of indices in $I_D\left(j\right)$ that are not multiple of 3:
\[A\left(j\right)=\left\{j\in I_D\left(j\right)\ |\ j\mathrm{\ mod\ }3\neq 0\right\}.\] 
Let also $B\left(j\right)$ be the indices among them corresponding to composite numbers:
\[B\left(j\right)=A\left(j\right)\cap W.\] 
\textbf{\underbar{}}

\noindent \textbf{\underbar{Property 2.1}}: A remarkable index cannot be of the form $3n+1$.

\noindent More precisely, $k_j\left(j+1\right)\mathrm{\ mod\ }3=\left\{ \begin{array}{c}
0 \\ 
2 \end{array}
\ \  \begin{array}{c}
\mathrm{if}\ j=0\ \mathrm{mod\ }3; \\ 
\mathrm{otherwise}.\ \ \ \ \ \ \ \  \end{array}
\right.$

\noindent \textit{\underbar{Proof}}: We have $k_j\left(j+1\right)=2j^2+6j+3=-j^2\mathrm{\ mod\ }3$, which yields the result.

\noindent \textbf{\underbar{}}

\noindent \textbf{\underbar{Corollary 2.1}}: $\left|A\left(j\right)\right|$ is always even.

\noindent \textit{\underbar{Proof}}: We deduce from property 2.1 that $3n+1\in A\left(j\right)$ if and only if $3n+2\in A\left(j\right)$. There are therefore as many numbers congruent to 1 as 2 modulo 3 in the set $A\left(j\right)$. Thus, its size must be even.

\noindent 
\subsection{Prime number indices}

\noindent 

\noindent We adopt the usual notation $ \pi\left(x\right)$ for the number of primes not greater than \textit{x}. We will also note $ \pi\left(j\right)$ for the number of indices corresponding primes not greater than ${\left(2j+5\right)}^2$, i.e. the elements in $A\left(j\right)$. Thus, by definition: \\
$ \pi{'} \left(j\right) = \pi \left({\left(2j+5\right)}^2\right)-2$
as 2 and 3 must be removed. \\ 

\noindent \textbf{\underbar{Property 2.2}:} We have the following equality:
\[\pi{'}\left(j\right)=\left|A\left(j\right)\right|-\left|B\left(j\right)\right|.\]
\textit{\underbar{Proof}}: As every odd number is either prime or composite, it is immediate that $A\left(j\right)\backslash B\left(j\right)$ is the set of prime number indices in $A\left(j\right)$, and the result follows.

\noindent 
\section{Counting $\boldsymbol{A}$($\boldsymbol{j}$) and $\boldsymbol{B}$($\boldsymbol{j}$) }

\noindent We will now proceed to counting the sets defined in the previous section.

\noindent 
\subsection{Counting $\boldsymbol{A}$($\boldsymbol{j}$)}

\noindent We start with a helpful lemma:

\noindent \textbf{\underbar{Lemma 3.1}}: Let $n\in {\mathbb{N}}^*$ and $X$ be set of consecutive integers (or \underbar{range}), the size of which is a multiple of $n$. Then:
\[\forall R\subset \llbracket 0,n-1\rrbracket \ \left|\left\{x\in X|\left(x\ \mathrm{mod}\ n\right)\in R\right\}\right|=\frac{\left|R\right|.\left|X\right|}{n}.\] 
\textit{\underbar{Proof}}: We take $n\in {\mathbb{N}}^*$ and $R\subset \llbracket 0,n-1\rrbracket $ fixed, and we proceed by induction on $\left|X\right|$.

\noindent The result is obviously true for $X=\emptyset $. Suppose now that it is also true for any range $X$ of size $an$ ($a\ge 0$). Let $X$ be a range of size $\left(a+1\right)n$, and $x_0$ be its smallest element. Let us then define $X_0$ and $X{'}$ by:
\[X_0=\llbracket x_0,x_0+n-1\rrbracket ,\] 
\[X{'}=X\backslash X_0.\] 
It is straightforward that $\left|X_0\right|=n$ and $\left|X{'}\right|=an$, therefore it follows from the induction hypothesis that,
\[\left|\left\{x\in X|\left(x\ \mathrm{mod}\ n\right)\in R\right\}\right|=\left|\left\{x\in X_0|\left(x\ \mathrm{mod}\ n\right)\in R\right\}\right|+\frac{\left|R\right|.\left|X'\right|}{n}\] 
But $X_0$ is a range of $n$ consecutive integers, thus each congruence class appears exactly once, from which we conclude that $\left|\left\{x\in X_0|\left(x\ \mathrm{mod}\ n\right)\in R\right\}\right|=\left|R\right|$. This proves the lemma.

\noindent 

\noindent For a number $x\in \mathbb{R}$, we will note $\left\lfloor x\right\rfloor $ its integral part, which is also the number of positive integers not greater than $x$.

\noindent \textbf{\underbar{Property 3.1}}: For all $j\in \mathbb{N}$, we have:
\[\left|A\left(j\right)\right|=\left\lfloor \frac{2}{3}\left(k_{j+1}\left(j+2\right)+1\right)\right\rfloor \] 
\textbf{\underbar{}}

\noindent \textit{\underbar{Proof}}:\underbar{}

\noindent On one hand, if $k_{j+1}\left(j+2\right)\ \mathrm{mod\ }3=2$, the result follows directly from the previous lemma.

\noindent On the other hand, we know that otherwise $k_{j+1}\left(j+2\right)$ has to be a multiple of 3, therefore $k_{j+1}\left(j+2\right)\notin A\left(j\right)$ and the same lemma applied to $\llbracket 0,k_{j+1}\left(j+2\right)-1\rrbracket $ yields:
\[\left|A\left(j\right)\right|=\frac{2}{3}\left(k_{j+1}\left(j+2\right)\right).\] 
We deduce that in all cases:
\[\left|A\left(j\right)\right|=\left\lfloor \frac{2}{3}\left(k_{j+1}\left(j+2\right)+1\right)\right\rfloor .\]

\noindent \textbf{\underbar{Corollary 3.1}}: We have the following asymptotic expansion:
\[\left|A\left(j\right)\right|=\frac{4}{3}j^2+\frac{20}{3}j+O\left(1\right).\] 
In particular:
\[\left|A\left(j\right)\right|\leadstoext \frac{4}{3}j^2.\]

\noindent 
\subsection{Counting $\boldsymbol{B}$($\boldsymbol{j}$)}

\noindent 
\subsubsection{The inclusion-exclusion principle}

\noindent We remind that $W$ is the union of the sets $W_j$ of indices corresponding to composite odd multiples of $2j+3$. In particular:
\[B_j=(W\cap I_D\left(j\right))\backslash \left(W_0\cap I_D\left(j\right)\right).\] 
Let $x\in W\cap I_D\left(j\right)$. Thus, $2x+3\le {\left(2j+5\right)}^2$, and $2x+3$ is a composite number, so it admits at least one odd prime factor not greater than $2j+5$. There are $\pi{''}\left(j\right)+1=\pi \left(2j+5\right)-1$ such primes (2 is excluded but not 3). We deduce that:
\[W\cap I_D\left(j\right)=\bigcup^{{\pi }^{''}\left(j\right)}_{k=0}{W_{q_k}\cap I_D\left(j\right)}.\] 
Furthermore, we may involve more prime numbers without changing the result. We deduce that more generally:
\[\forall N\ge \pi{''}\left(j\right)\ \ W\cap I_D\left(j\right)=\bigcup^N_{k=0}{W_{q_k}\cap I_D\left(j\right)}.\]

\noindent The inclusion-exclusion principle implies that the size of $B_j$ verifies:
\[ \forall N\ge \pi{''}\left(j\right)\ \ \left|B_j\right|=\sum_{ K\subset \llbracket 0,N\rrbracket ,\ K\neq \emptyset ,\ K\neq \left\{0\right\}}
{{\left(-1\right)}^{\left|K\right|-1}\left|\left(\bigcap_{k\in K}{W_{q_k}}\right)\cap I_D\left(j\right)\right|}.\#\left(1\right)\] 
\textbf{\underbar{Definition 3.2.1}}\underbar{:} Let $k,l$ be the indices of two odd numbers, we note $k*l$ for the index of their product. We verify that:
\[k*l=2kl+3\left(k+l+1\right)\] 
We know that the product of integers is associative and commutative. This implies that $*$ has these properties too.

\noindent For a set of any integers $K=\left\{i_1\dots i_n\right\}$ we note $q_K=q_{i_1}*\dots *q_{i_n}$. This definition is non-ambiguous because $*$ is associative. Furthermore, we can rewrite (1) as:
\[ \forall N\ge \pi{''}\left(j\right)\ \ \left|B_j\right|=\sum_{K\subset \llbracket 0,N\rrbracket  ,\ K\neq \emptyset ,\ K\neq \left\{0\right\}
}{{\left(-1\right)}^{\left|K\right|-1}\left|W_{q_K}\cap I_D\left(j\right)\right|}\#\left(2\right)\] 
Taking $N$ sufficiently large, all indices of \underbar{square-free} odd numbers between $5$ and ${\left(2j+5\right)}^2$ appear in the sum above, while all the indices greater to $k_{j+1}\left(j+2\right)$ yield no contribution to the sum (because if $q>k_{j+1}\left(j+2\right)$, obviously $\left|W_q\cap I_D\left(j\right)\right|=0$). We deduce a third version of (2) using the function of M\"{o}bius (see [2]), as we observe that ${\left(-1\right)}^{\left|K\right|}=\mu \left(2q_K+3\right)$:
\[ \begin{array}{c}
\left|B_j\right|=-\sum^{k_{j+1}\left(j+2\right)}_{k=1}{\mu \left(2k+3\right)\left|W_k\cap I_D\left(j\right)\right|}\ \#\left(3\right) \end{array}
\] 

\subsubsection{Calculation of the cardinality of  ${\boldsymbol{W}}_{\boldsymbol{k}}\boldsymbol{\cap }{\boldsymbol{I}}_{\boldsymbol{D}}$($\boldsymbol{j}$)}

\noindent 

\noindent Once again we start with a useful counting lemma:

\noindent \textbf{\underbar{Lemma 3.2}}: For all $n,m\in {\mathbb{N}}^*$ the number of multiples of $n$ between $1$ and $n$ is equal to:
\[\left|\llbracket 1,m\rrbracket \cap n\mathbb{Z}\right|=\left\lfloor \frac{m}{n}\right\rfloor .\] 
\textit{\underbar{Proof}}: We take $n\in {\mathbb{N}}^*$, and we proceed by induction on $m$. For $m=1$ the result is trivially true. Suppose that for a given value of $m\in {\mathbb{N}}^*$, we have $\left|\llbracket 1,m\rrbracket \cap n\mathbb{Z}\right|=\left\lfloor \frac{m}{n}\right\rfloor $. Then, if $m+1$ is not a multiple of $n$, we have:
\[\left|\llbracket 1,m+1\rrbracket \cap n\mathbb{Z}\right|=\left|\llbracket 1,m\rrbracket \cap n\mathbb{Z}\right|=\left\lfloor \frac{m}{n}\right\rfloor =\left\lfloor \frac{m+1}{n}\right\rfloor .\] 
On the other hand, if $m+1$ is a multiple of $n$, we deduce:
\[\left|\llbracket 1,m+1\rrbracket \cap n\mathbb{Z}\right|=1+\left|\llbracket 1,m\rrbracket \cap n\mathbb{Z}\right|=1+\left\lfloor \frac{m}{n}\right\rfloor =\left\lfloor \frac{m+1}{n}\right\rfloor .\] 
In all cases we manage to show that the property is inductive, which yields the result.

\noindent  \textbf{\underbar{Property 3.2.2:}} Let $k$ be an index. We have the following:
\[\left|W_k\cap I_D\left(j\right)\right|=\left\lfloor \frac{k_{j+1}\left(j+2\right)-k}{2k+3}\right\rfloor \] 
\textit{\underbar{Proof}}: Let $x\in \mathbb{N}$. Then $x\in W_k\Leftrightarrow \exists n\in {\mathbb{N}}^*\ \ x=\left(2k+3\right)n+k\Leftrightarrow x>k\ \mathrm{and\ }\left(2k+3\right)|\left(x-k\right)$.

\noindent Therefore, $x\in W_k\cap I_D\left(j\right)\Leftrightarrow 1\le x-k\le k_{j+1}\left(j+2\right)-k\ \mathrm{and\ }\left(2k+3\right)|\left(x-k\right)$

\noindent So $\left|W_k\cap I_D\left(j\right)\right|=\left|\left(W_k\cap I_D\left(j\right)\right)-k\right|=\left|\llbracket 1,k_{j+1}\left(j+2\right)-k\rrbracket \cap \left(2k+3\right)\mathbb{Z}\right|$.

\noindent The result then follows from lemma 3.2.

\noindent 

\noindent \textbf{\underbar{Corollary 3.2.2:}} If $n=2k+3$, we may also write:
\[\left|W_k\cap I_D\left(j\right)\right|=\left\lfloor \frac{{\left(2j+5\right)}^2-n}{2n}\right\rfloor \] 
\textit{\underbar{Proof}}: Indeed, we recall that by definition of the remarkable index, ${\left(2j+5\right)}^2=2k_{j+1}\left(j+2\right)+3$. Thus the result follows.

\noindent 

\noindent \textbf{\underbar{Definition~3.2.2:}} We now note $p_K=\prod_{n\in K}{p_n}$ the odd number indexed by $q_K$.

\noindent 

\noindent The property 3.2.3 gives two expressions for $\left|B_j\right|$:

\noindent \textbf{\underbar{Property 3.2.3:}}
\[\left|B_j\right|=\sum_{ \begin{array}{c}
K\subset \llbracket 0,\pi{''}\left(j\right)\rrbracket  \\ 
K\neq \emptyset ,\ K\neq \left\{0\right\} \end{array}
}{{\left(-1\right)}^{\left|K\right|-1}\left\lfloor \frac{{\left(2j+5\right)}^2-p_K}{2p_K}\right\rfloor }\] 

\[\left|B_j\right|=-\sum^{k_{j+1}(j+2)}_{k=1}{\mu \left(2k+3\right)\left\lfloor \frac{k_{j+1}\left(j+2\right)-k}{2k+3}\right\rfloor }\]

\noindent \textit{\underbar{Proof}}: It stems from $\left(2\right)$ and $\left(3\right)$, to which we apply property 3.2.2 and its corollary.

\noindent \textit{\underbar{Remark:}} In the first expression, we may group the terms by the size of $K$, which leads to a further expression with an alternating sum:
\[\left|B_j\right|=\left(\sum^{\pi{''}\left(j\right)}_{k=1}{\left\lfloor \frac{{\left(2j+5\right)}^2-p_k}{2p_k}\right\rfloor }\right)+\sum^{\pi{''}\left(j\right)}_{n=2}{{\left(-1\right)}^{n-1}\left(\sum_{ \begin{array}{c}
K\subset \llbracket 0,\pi{''}\left(j\right)\rrbracket  \\ 
\left|K\right|=n \end{array}
}{\left\lfloor \frac{{\left(2j+5\right)}^2-p_K}{2p_K}\right\rfloor }\right)}\] 

\subsubsection{Asymptotic expansion of {\textbar}${\boldsymbol{B}}_{\boldsymbol{j}}${\textbar}}

\noindent The prime number theorem [3], demonstrated independently by Hadamard and de la Vall\'{e}e Poussin in 1896, is an important result on the asymptotic expansion of the number of prime numbers. It states that for $x\to +\infty $: 
\[\pi \left(x\right)\sim \frac{x}{{\mathrm{ln} \left(x\right)\ }}\] 
\textbf{\underbar{Property 3.2.4}}: We have the following asymptotic expansion:
\[\left|B_j\right|=\frac{4}{3}j^2-\frac{2j^2}{{\mathrm{ln} \left(j\right)\ }}+o\left(\frac{j^2}{{\mathrm{ln} \left(j\right)\ }}\right)\] 
\textit{\underbar{Proof}}: Property 2.2 gives $\left|B_j\right|=\left|A_j\right|-\pi{'}\left(j\right)$. Corollary 3.1 gives a very precise asymptotic expansion of $\left|A_j\right|$. From the prime number theorem, we also deduce that $\pi{'}\left(x\right)=\pi \left({\left(2j+5\right)}^2\right)-2\leadstoext \frac{2j^2}{{\mathrm{ln} \left(j\right)\ }}$. Thus we conclude.

\noindent 

\noindent \underbar{Remark:} Even with the known refinements of the prime number theorem, it is not possible to improve the result in $O\left(j\right)$, let alone $O\left(1\right)$.

\noindent 
\subsection{A special alternate series equivalent to {\textbar}${\boldsymbol{B}}_{\boldsymbol{j}}${\textbar} }

\noindent In this section we focus on another equivalent of $\left|B_j\right|$.

\noindent A naive manipulation of the formula derived from property 3.2.3:
\[\left|B_j\right|=\left(\sum^{\pi{''}\left(j\right)}_{k=1}{\left\lfloor \frac{{\left(2j+5\right)}^2-p_k}{2p_k}\right\rfloor }\right)+\sum^{\pi{''}\left(j\right)}_{n=2}{{\left(-1\right)}^{n-1}\left(\sum_{ \begin{array}{c}
K\subset \llbracket 0,\pi{''}\left(j\right)\rrbracket  \\ 
\left|K\right|=n \end{array}
}{\left\lfloor \frac{{\left(2j+5\right)}^2-p_K}{2p_K}\right\rfloor }\right)}\] 
would consist in summing asymptotic equivalents of each term of the sum, which would lead to an expression without integral parts:
\[\left|B_j\right|\leadstoext 2\left(\sum^{\pi{''}\left(j\right)}_{k=1}{\frac{1}{p_k}}+\sum^{\pi{''}\left(j\right)}_{n=2}{{\left(-1\right)}^{n-1}\sum_{ \begin{array}{c}
K\subset \llbracket 0,\pi{''}\left(j\right)\rrbracket  \\ 
\left|K\right|=n \end{array}
}{\frac{1}{p_K}}}\right)j^2.\] 
However, we must be careful that, as the number of terms in the sum is not bounded, this approach is not mathematically valid. We may however show using Eulerian products that the result \underbar{is} correct.

\noindent \textbf{\underbar{Property 3.3.1}}: The coefficient $c_j=\sum^{\pi{''}\left(j\right)}_{k=1}{\frac{1}{p_k}}+\sum^{\pi{''}\left(j\right)}_{n=2}{{\left(-1\right)}^{n-1}\sum_{ \begin{array}{c}
K\subset \llbracket 0,\pi{''}\left(j\right)\rrbracket  \\ 
\left|K\right|=n \end{array}
}{\frac{1}{p_K}}}$ converges to ${2}/{3}$ when $j\to +\infty $. 

\noindent \textit{\underbar{Proof}}: Euler (see [4]) proved the divergence of the series of the reciprocals of the primes: 
\[\sum{\frac{1}{p_i}}=+\infty \] 
As $-{\mathrm{ln} \left(1-\frac{1}{p_i}\right)\ }\leadstoext \frac{1}{p_i}$, the limit comparison test shows that:
\[\sum{{\mathrm{ln} \left(1-\frac{1}{p_i}\right)\ }}=-\infty \] 
Therefore, using the exponential function:
\[\prod{\left(1-\frac{1}{p_i}\right)}=0\] 
By developing the finite version of the product above, our alternating sum almost appears:
\[\prod^n_{i=0}{\left(1-\frac{1}{p_i}\right)}=1-\sum^n_{i=0}{\frac{1}{p_i}}+\sum_{0\le i<j\le n}{\frac{1}{p_ip_j}}-\dots +\frac{{\left(-1\right)}^n}{p_1\dots p_n}\] 
Remember that $c_j$ is equal to:
\[\sum^n_{i=1}{\frac{1}{p_i}}-\sum_{0\le i<j\le n}{\frac{1}{p_ip_j}}-\dots +\frac{{\left(-1\right)}^{n+1}}{p_0\dots p_n}\] 
The limit of the former being zero, the latter therefore converges to the sum of the terms removed, i.e. ${\mathrm{lim}}$ $c_j\ =1-{1}/{3}={2}/{3}$.

\noindent \textbf{\underbar{Corollary 3.3.1}}: We deduce that:
\[\left|B_j\right|\leadstoext 2\left(\sum^{\pi{''}\left(j\right)}_{k=1}{\frac{1}{p_k}}+\sum^{\pi{''}\left(j\right)}_{n=2}{{\left(-1\right)}^{n-1}\sum_{ \begin{array}{c}
K\subset \llbracket 0,\pi{''}\left(j\right)\rrbracket  \\ 
\left|K\right|=n \end{array}
}{\frac{1}{p_K}}}\right)j^2.\] 
Furthermore:
\[c_j=\sum^{\pi{''}\left(j\right)}_{k=1}{\frac{1}{p_k}}+\sum^{\pi{''}\left(j\right)}_{n=2}{{\left(-1\right)}^{n-1}\sum_{ \begin{array}{c}
K\subset \llbracket 0,\pi{''}\left(j\right)\rrbracket  \\ 
\left|K\right|=n \end{array}
}{\frac{1}{p_K}}}=\frac{2}{3}\left(\sum^{\pi{''}\left(j\right)}_{n=1}{{\left(-1\right)}^{n-1}\sum_{ \begin{array}{c}
K\subset \llbracket 1,\pi{''}\left(j\right)\rrbracket  \\ 
\left|K\right|=n \end{array}
}{\frac{1}{p_K}}}\right)\] 
\textit{\underbar{Proof}}: Property 3.3.1 shows that $c_j\leadstoext {2}/{3}$, and property 3.2.4 implies $\left|B_j\right|\leadstoext {4j^2}/{3}$. Thus, $\left|B_j\right|\leadstoext 2c_j.j^2$ which proves the first part of the corollary. The second part is an alternative expression of $c_j$ obtained by isolating every term containing 3 in the sum: \\
\[\sum^{\pi{''}\left(j\right)}_{k=1}{\frac{1}{p_k}}+\sum^{\pi{''}\left(j\right)}_{n=2}{{\left(-1\right)}^{n-1}\sum_{ \begin{array}{c}
K\subset \llbracket 0,\pi{''}\left(j\right)\rrbracket  \\ 
\left|K\right|=n \end{array}
}{\frac{1}{p_K}}}=\] \\
\[\sum^{\pi{''}\left(j\right)}_{k=1}{\frac{1}{p_k}}+\sum^{\pi{''}\left(j\right)}_{n=2}{{\left(-1\right)}^{n-1}\sum_{ \begin{array}{c}
K\subset \llbracket 1,\pi{''}\left(j\right)\rrbracket  \\ 
\left|K\right|=n \end{array}
}{\frac{1}{p_K}}}+\sum^{\pi{''}\left(j\right)}_{n=2}{{\left(-1\right)}^{n-1}\sum_{ \begin{array}{c}
K\subset \llbracket 0,\pi{''}\left(j\right)\rrbracket  \\ 
\left|K\right|=n \\ 
0\in K \end{array}
}{\frac{1}{p_K}}}=\] \\
\[\left(1-\frac{1}{p_0}\right).\left(\sum^{\pi{''}\left(j\right)}_{n=1}{{\left(-1\right)}^{n-1}\sum_{ \begin{array}{c}
K\subset \llbracket 1,\pi{''}\left(j\right)\rrbracket  \\ 
\left|K\right|=n \end{array}
}{\frac{1}{p_K}}}\right).\]

The last equality is obtained writing $p_K=p_0.p_{K{'}}$ for any $K$ containing $0$.

\noindent 
\subsection{Another interesting limit  }

\noindent We get a similar result for the M\"{o}bius version of $\left|B_j\right|$ in 3.2.3:

\noindent \textbf{\underbar{Property 3.3.2}}:
\[\sum^{+\infty }_{n=2}{\frac{-\mu \left(2n+1\right)}{2n+1}}=\frac{2}{3}.\] 
\textit{\underbar{Proof}}: The prime number theorem is the equivalent to the following (see [5]):
\[\sum^{+\infty }_{n=1}{\frac{\mu \left(n\right)}{n}}=0.\] 
The even terms $\frac{\mu \left(2n\right)}{2n}$ are undesirable, but we note that $\mu \left(2n\right)$ is non-zero only if $n$ is odd, in which case $\mu \left(2n\right)=-\mu \left(n\right)$. Let $S_N=\sum^N_{n=1}{\frac{\mu \left(n\right)}{n}}$ and $T_N=\sum^{N-1}_{n=0}{\frac{\mu \left(2n+1\right)}{2n+1}}$. We have:
\[S_{4N}=\sum_{ \begin{array}{c}
n\le 4N \\ 
n\ \mathrm{odd} \end{array}
}{\frac{\mu \left(n\right)}{n}}+\sum_{ \begin{array}{c}
n\le 4N \\ 
n\ \mathrm{even} \end{array}
}{\frac{\mu \left(n\right)}{n}}=T_{2N}-\sum_{ \begin{array}{c}
n\le 2N \\ 
n\ \mathrm{odd} \end{array}
}{\frac{\mu \left(n\right)}{2n}}=T_{2N}-\frac{T_N}{2}.\] 
It follows that if $a$ a is a cluster point of $\left(T_{2N}\right)$, $2a$ is necessarily a cluster point of $\left(T_N\right)$, as $S_{4N}$ converges to $0$. But the difference $T_{N+1}-T_N$ converges also to $0$, which yields that the cluster points of $\left(T_{2N}\right)$ and $\left(T_N\right)$ are the same, and that they form a range $A$, with the property $a\in A\Rightarrow 2a\in A$.

\noindent All that remains is to prove that $\left(T_N\right)$ is bounded, which will yield that $A$ is necessarily equal to $\left\{0\right\}$, and that $\left(T_N\right)$ converges to 0.

\noindent To this end, we need to extend property 3.2.3: 
\[\left|B_j\right|=\sum^{k_{j+1}\left(j+2\right)}_{k=1}{\mu \left(2k+3\right)\left\lfloor \frac{k_{j+1}\left(j+2\right)-k}{2k+3}\right\rfloor }.\] 
Indeed, we notice that every odd number greater than one is a multiple of at least one odd prime, which yields:
\[\left\{3,5,\dots ,2n+1\right\}=\bigcup^{+\infty }_{k=0}{p_k\left\{1,3,\dots ,2\left\lfloor \frac{2n+1}{2p_k}-\frac{1}{2}\right\rfloor +1\ \right\}}.\] 
Thus, using the inclusion-exclusion formula:
\[\forall n\in \mathbb{N}\mathrm{\ }\mathrm{\ }n=\sum^n_{k=1}{-\mu \left(2k+1\right)\left\lfloor \frac{n+k+1}{2k+1}\right\rfloor }\] 
which is equivalent to:
\[\forall n\in \mathbb{N}\mathrm{\ }\mathrm{\ }1=\sum^n_{k=0}{\mu \left(2k+1\right)\left\lfloor \frac{n+k+1}{2k+1}\right\rfloor }.\] 
The following inequality follows from neglecting the integer parts (as the first and the last term are already integers, the error is at most $n-1$):
\[\forall n\in {\mathbb{N}}^*\ \ -n+2\le \sum^n_{k=0}{\mu \left(2k+1\right)\frac{n+k+1}{2k+1}}\le n\] 
Finally, as $\frac{k+1}{2k+1}<1$ for all $k>0$:
\[\forall n\in {\mathbb{N}}^*\ \ \left|\sum^n_{k=0}{\frac{\mu \left(2k+1\right)}{2k+1}}\right|<2.\] 
This proves the boundedness of $T_N$.

\noindent Eventually, from $\sum^{+\infty }_{n=0}{\frac{\mu \left(2n+1\right)}{2n+1}}=0$ it is easy to deduce that:
\[\sum^{+\infty }_{n=2}{\frac{-\mu \left(2n+1\right)}{2n+1}}=1-\frac{1}{3}=\frac{2}{3}.\] 
\underbar{Remark}: From the relationship $S_{4N}=T_{2N}-\frac{T_N}{2}$, we can conversely deduce that if $T_N$ converges, so does $S_N$, which proves the convergence of $T_N$ is equivalent to the prime number theorem.

\noindent 
\subsection{Convergence rate comparison }

\noindent The sum in property 3.3.2 contains fewer terms than that of property 3.3.1, and terms are summed in a different order. These two convergence results are therefore not equivalent. In this last part, we will empirically compare the behavior of these two sums with that of $\left|B_j\right|$ and its asymptotic expansion given directly by the prime number theorem (property 3.2.4).

\noindent To make everything comparable, we will normalize all these quantities so that they represent approximate proportions of the composite numbers among odd numbers non-multiple of 3.

\noindent We define 4 sequences:

\begin{enumerate}
\item  $p_j=\left|B_j\right|/\left|A_j\right|$ the \underbar{exact proportion};

\item  $a_j=1-\frac{3}{2{\mathrm{ln} \left(j\right)\ }}$ the \underbar{approximate proportion} to order 1 (or \underbar{Hadamard approximation});

\item  \[e_j=\sum^{\pi{''}\left(j\right)}_{n=1}{{\left(-1\right)}^{n-1}\sum_{ \begin{array}{c}
K\subset \llbracket 1,\pi{''}\left(j\right)\rrbracket  \\ 
\left|K\right|=n \end{array}
}{\frac{1}{p_K}}}=1-\prod^{\pi{''}\left(j\right)}_{k=1}{\left(1-\frac{1}{p_k}\right)}\] the \underbar{Euler approximation}.

\item  $m_j=\frac{3}{2}\sum^{k_{j+1}\left(j+2\right)}_{k=1}{\frac{-\mu \left(2k+3\right)}{2k+3}}$ the \underbar{M\"{o}bius approximation}.
\end{enumerate}

\noindent In the graph below, we set the squares${\left(2j+5\right)}^2$ on the x-axis (with a logarithmic scale), and on the corresponding proportions $p_j,\ a_j,\ e_j$ et $m_j$ on the y-axis:

\noindent 
\includegraphics[width=\linewidth]{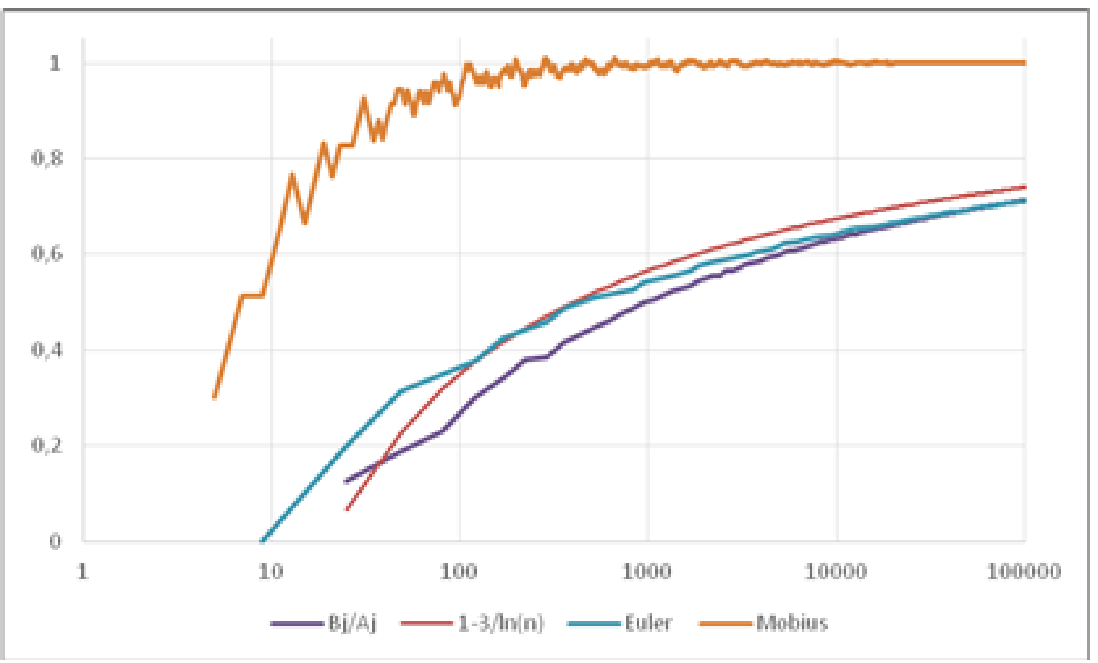} 
\noindent 

\noindent The result suggest that the Euler approximation is the best of the three, whereas the M\"{o}bius approximation converges much faster to 1 than the others, which also makes it a poorer approximation. We also observe that it is the least monotonous and, in terms of complexity, the one that requires most operations (its expression requires to compute the M\"{o}bius function of every number between 1 and ${\left(2j+5\right)}^2$, i.e. a full factorization for all square-free numbers).

\noindent 
\part*{Conclusion}
\addcontentsline{toc}{part}{Conclusion}
\noindent With some simple counting arguments applied to finite sets of composite odd numbers (via their indices), two harmonic sums appeared naturally. We proved their convergence which is also the illustration of the fact there are asymptotically almost as many composite odd numbers as odd numbers -- the sums were shown to be approximations of the exact ratio. One of these convergences was also proven to be equivalent to the prime number theorem. An empirical discussion of the quality of these approximations suggests that there is however a profound difference between the two.

\noindent 
\part*{References}
\addcontentsline{toc}{part}{References}
\noindent [1] WOLF Marc, WOLF Fran\c{c}ois, Representation theorem of composite odd numbers indices, SCIREA Journal of Mathematics. Vol. 3, No. 3, 2018, pp. 106 - 117.

\noindent [2] G. Tenenbaum et M. Mend\`{e}s-France, Les nombres premiers, entre l'ordre et le chaos, p. 25

\noindent [3] G. H. Hardy et E. M. Wright, An Introduction to the Theory of Numbers, 4e \'{e}d., p. 10.

\noindent [4] G. Tenenbaum et M. Mend\`{e}s-France, Les nombres premiers, entre l'ordre et le chaos, p. 23

\noindent [5] G. Tenenbaum et M. Mend\`{e}s-France, Les nombres premiers, entre l'ordre et le chaos, p. 126

\end{document}